\documentclass[11pt,twoside,a4paper]{article}

\usepackage{amsmath}
\usepackage{amsthm}
\usepackage{changepage}
\usepackage{geometry}
\usepackage{txfonts}
\usepackage{cite} 
\usepackage{listings}
\usepackage{hyperref}
\usepackage{amssymb}
\usepackage{mathrsfs}
\author{}
\hypersetup{
	colorlinks=true,
	linkcolor=blue,
	anchorcolor=blue,
	citecolor=blue}
\title{ \textbf{Relation between the eventual continuity and the e-property for Markov-Feller semigroups}}
\newtheorem{theorem}{Theorem}[section]
\newtheorem{lemma}[theorem]{Lemma}
\newtheorem{corollary}[theorem]{Corollary}
\newtheorem{proposition}[theorem]{Proposition}
\numberwithin{equation}{section}
\theoremstyle{definition}
\newtheorem{remark}[theorem]{Remark}
\newtheorem{example}[theorem]{Example}
\newtheorem{definition}[theorem]{Definition}
\newcommand{\runum}[1]{\romannumeral #1}
\newenvironment{prof}[1][Proof]{\noindent\textit{#1}\quad }

\begin{document}

 {\author{ \textbf{Yong Liu,\   Ziyu Liu} \\ 
\small{\em     LMAM, School of Mathematical Science, Peking University, 100871, Beijing, China}\\
\small{\em liuyong@math.pku.edu.cn (Y. Liu);
\ ziyu@pku.edu.cn(Z. Liu)}
}}

\date{February 17, 2023}
	\maketitle

	\makeatother
	
	\begin{abstract}
	We investigate the relation between the e-property and the eventual continuity, or called the asymptotic equicontinuity, which is a generalization of the e-property.  We prove that, for any  discrete-time or strongly continuous continuous-time eventually continuous  Markov-Feller semigroup with an ergodic measure,  if the interior of the support of the ergodic measure is nonempty, then the e-property is satisfied on the interior of the support. In particular, it implies that, restricted on the support of each ergodic measure, the e-property and the eventual continuity are equivalent for the discrete-time and the strongly  continuous continuous-time Markov-Feller semigroups.
    \end{abstract}
\renewcommand{\thefootnote}{}
\footnotetext{Y. Liu is supported by CNNSF (No. 11731009,  No. 12231002) and Center for Statistical Science, PKU. } 

     \section{Introduction}\label{Sec 1}
      In this paper, we are mainly concerned with the  relation between the regularities  of the Markov-Feller semigroups. In 2006, Lasota and Szarek set forth the concept of \emph{e-property} and developed the lower bound technique to formulate a criterion for the existence of an invariant measure\cite{S2000,LY1994}. Since then, the e-property has turned out to be a valuable tool in proving the existence of invariant measures and the ergodicity, and abundant results have been obtained for both the Markov semigroups and their Ces\`aro averages\cite{SW2012,KSS2012,HSZ2017,WW2018,BKS2014}.  \par 

     This paper is mainly motivated by the analysis of the non-equicontinuous Markov semigroups, that is, Markov semigroups which do not satisfy the e-property. This kind of Markov semigroups can be usually discovered in many stochastic dynamical systems.  For example, non-equicontinuous Markov-Feller semigroups, generated by deterministic dynamical systems, have been given in \cite{GL2015,HSZ2017,J20132}. To handle the ergodicity of the non-equicontinuous Markov semigroups, Czapla introduced a generalization of the e-property for the Markov-Feller semigroups in \cite{C2017}, the \emph{eventual e-property}, which is weaker than the e-property. Indeed, we  show that these two properties are equivalent for the Markov-Feller semigroups (see Proposition \ref{Prop 1}, \ref{Prop 2}).   Another even weaker condition, the notion of  the \emph{asymptotic equicontinuity condition} was introduced by Jaroszewska in \cite{J2013}, or called the \emph{eventual continuity} by Gong and Liu in \cite{GL2015}.  In fact, these two notions are formulated almost simultaneously, and mathematically equivalent. We adopt the name of the eventual continuity in this paper.   The  eventual continuity depicts the feature that a uniquely ergodic semigroup may behave sensitively in initial data, and is evidently weaker than the  e-property.  The  eventual continuity is a necessary condition for both the e-property and the \emph{asymptotic stability}, and is a reasonable tool to investigate the ergodicity of the non-equicontinuous Markov-Feller semigroups, see  \cite{J2013,GL2015,GLLL2023}.\par

    One natural question arises: when does a general Markov-Feller semigroup satisfy the e-property?   In \cite{HSZ2017}, Hile, Szarek and Ziemla\'nska showed that any asymptotically stable discrete-time  Markov-Feller semigroup with an invariant measure such that the interior of its support is nonempty satisfied the e-property. In this paper, instead of the asymptotic stability, we use the eventual continuity to provide 
    a more general criterion for the e-property.  Specifically, when a discrete-time or strongly continuous continuous-time eventually continuous Markov-Feller semigroup has an ergodic measure  such that the interior of its support is nonempty, then it satisfies the e-property on the interior of the support of the ergodic measure.    This result implies that the e-property and the eventual continuity are equivalent, restricted on the support of each ergodic measure of any discrete-time or  strongly continuous continuous-time Markov-Feller semigroup.   \par

    The paper is organized as follows. Section \ref{Sec 2} contains some notions and definitions from the theory of Markov operators and the semigroup regularities.  The main results are provided in Section \ref{Sec 3}. Subsection \ref{Sec 3.1} deals with the relation between the e-property and the eventual e-property. Subsection \ref{Sec 3.2} is devoted for the relation between the e-property and the eventual continuity. Subsection \ref{Sec 3.3} provides some further conclusions and discussions  on the semigroup regularities.  The proofs of lemmas are placed in Section \ref{Sec 4}.

    \section{Preliminaries}\label{Sec 2}
    Let $(\mathcal{X},\rho)$ be a Polish space, i.e. a separable, complete metric space, and $\mathcal{B}(\mathcal{X})$ denote the $\sigma$-field of all its Borel subsets. We introduce the following notations: 
	
	$\begin{aligned}
        \mathcal{M}(\mathcal{X})&=\text{ the family of all finite Borel measures on } \mathcal{X},\\
	\mathcal{P}(\mathcal{X})&=\text{ the family of all probability measures on } \mathcal{X},\\
	B_b(\mathcal{X})&=\text{ the space of all bounded, Borel real-valued functions defined on }X,\\
	&\quad\;\text { endowed with the supremum norm: } \|f\|_{\infty}=\text{sup}_{x\in \mathcal{X}}\,|f(x)|,f\in B_b(\mathcal{X}),\\
	C_b(\mathcal{X})&=\text{ the subspace of }  B_b(\mathcal{X}) \text{ consisting of all bounded continuous functions},\\
	L_b(\mathcal{X})&=\text{ the subspace of } C_b(\mathcal{X}) \text{ consisting of all bounded Lipschitz functions},\\
	B(x,r)&=\{y\in \mathcal{X}:\rho(x,y)<r\} \text{ for } x\in \mathcal{X} \text{ and }r>0,\\
	\partial A,\overline{A},\text{Int}_{\mathcal{X}} (A)&=  \text{ the boundary, closure, interior of } A \text{ in } \mathcal{X}, \text{ respectively}, \\
	\text{supp }\mu&=\{x\in \mathcal{X}:\mu(B(x,\epsilon))>0 \text{ for every } \epsilon>0\}, \text{ for }\mu\in\mathcal{M}(\mathcal{X}), \\               &\quad\text{ i.e. the support of the measure } \mu ,\\
        T&=\text{ the index set, }\mathbb{R}_+=[0,\infty)\text{ or }\mathbb{N}_+=\{1,2,\dots,\}.
	\end{aligned}$	
   
   \vskip2mm
   For brevity, we use the notation $\langle f,\mu\rangle=\int_{\mathcal{X}}f(x)\mu(dx)$  for $f\in B_b(\mathcal{X})$ and $\mu\in\mathcal{M}(\mathcal{X})$.	\par

    An operator $P:\mathcal{M}(\mathcal{X})\rightarrow\mathcal{M}(\mathcal{X})$ is called a Markov operator on $\mathcal{X}$ if it satisfies that
    \begin{itemize}
        \item[$(\runum{1})$] (Positive linearity) $P(\lambda_1\mu_1+\lambda_2\mu_2)=\lambda_1P\mu_1+\lambda_2P\mu_2$ for $\lambda_1,\lambda_2\geq 0$, $\mu_1,\mu_2\in\mathcal{M}(\mathcal{X})$;
        \item[$(\runum{2})$]  (Preservation of the norm) $P\mu(\mathcal{X})=\mu(\mathcal{X})$ for $\mu\in\mathcal{M}(\mathcal{X})$.
    \end{itemize}
    A Markov operator $P$ is called regular if there exists a linear operator $P^*:B_b(\mathcal{X})\rightarrow B_b(\mathcal{X})$ such that
    \begin{equation*}
        \langle f, P\mu\rangle = \langle P^*f, \mu\rangle\quad\text{for all }f\in B_b(\mathcal{X}),\;\mu\in\mathcal{M}(\mathcal{X}).
    \end{equation*}
    For ease of notation, we simply rewrite $P^*$ as $P$. 
     A Markov operator $P$ is called a Markov-Feller operator if it is regular and $P$ leaves $C_b(\mathcal{X})$ invariant, i.e., $P(C_b(\mathcal{X}))\subset C_b(\mathcal{X})$.   A Markov semigroup  $\{P_t\}_{t\in T}$ on  $\mathcal{X}$ is a semigroup of Markov operators on $\mathcal{M}(\mathcal{X})$. A Markov semigroup $\{P_t\}_{t\geq 0}$ is called a Markov-Feller semigroup if $P_t$ is a Markov-Feller operator for all $t\geq 0$.  Recall that $\mu\in\mathcal{P}(\mathcal{X})$ is invariant for the semigroup $\{P_t\}_{t\in T}$ if $P_t\mu=\mu$ for all $t\in T$.  For $\mu\in\mathcal{P}(\mathcal{X}),$ define 
     \begin{equation*}
         Q_t\mu:=\frac{1}{t}\sum_{k=1}^{t}P_k\mu\quad\text{for }T=\mathbb{N}_+;\quad Q_t\mu:=\frac{1}{t}\int_{0}^{t}P_s\mu ds\quad\text{for }T=\mathbb{R}_+,
     \end{equation*}
     and denote $Q_t(x,\cdot)=Q_t\delta_x.$ \par

    Throughout this paper, we assume that $\{P_t\}_{t\in T}$ is a Markov-Feller semigroup.  Recall some types of regularities of the Markov semigroups. \par 
    
    \begin{definition}\label{def e-property}
     A Markov semigroup $\{P_t\}_{t\in T}$ satisfies the  e-property (see \cite{MT1993}, for example) at $z\in\mathcal{X}$, if for every $f\in L_b(\mathcal{X})$ 
	\begin{equation*}
	\limsup\limits_{x\rightarrow z}\sup\limits_{t\in T}|P_tf(x)-P_tf(z)|=0,
	\end{equation*}
    that is, $\forall\,\epsilon>0$, $\exists\,\delta>0$, such that $\forall\,x\in B(z,\delta)$, and $t\in T$, $|P_tf(x)-P_tf(z)|<\epsilon$.
    \end{definition}

     \begin{definition}\label{def eventual e-property}
    A Markov semigroup $\{P_t\}_{t\in T}$ satisfies the eventual e-property (see \cite{C2017}) at $z\in\mathcal{X}$, if for every $f\in L_b(\mathcal{X})$ 
		\begin{equation*}
		\limsup\limits_{(x,t)\rightarrow(z,\infty)}|P_tf(x)-P_tf(z)|=0,
		\end{equation*}
    that is, $\forall\,\epsilon>0$, $\exists\,\delta>0$, $t_0\in T$, such that $\forall\,x\in B(z,\delta)$, and $t\geq t_0$, $|P_tf(x)-P_tf(z)|<\epsilon$.
     \end{definition}
     \begin{remark}
        As far as we know, the notion of the eventual e-property was first formulated by Worm in \cite{W2010} as follows:  
       a Markov-Feller semigroup $\{P_t\}_{t\in T}$ satisfies the \emph{eventual e-property}  if there exists $t_0\in T$ such that $\{P_tf\}_{t\geq t_0}$ is equicontinuous for every $f\in L_b(\mathcal{X})$.  Worm's definition is slightly stronger than that in \cite{C2017}  (also see Definition \ref{def eventual e-property}).
     \end{remark}
     
    \begin{definition}\label{def eventual continuity}
         A Markov semigroup $\{P_t\}_{t\in T}$ satisfies the eventual continuity (see \cite{GL2015}) or the asymptotic equicontinuity (see \cite{J2013}) at $z\in\mathcal{X}$, if for every $f\in L_b(\mathcal{X})$ 
	\begin{equation*}
	\limsup\limits_{x\rightarrow z}\limsup\limits_{t\rightarrow \infty}|P_tf(x)-P_tf(z)|=0,
	\end{equation*}
	that is, $\forall\,\epsilon>0$, $\exists\,\delta>0$,  such that $\forall\,x\in B(z,\delta)$, $\exists\,t_x\in T$,   $\forall\,t\geq t_x$, $|P_tf(x)-P_tf(z)|<\epsilon$.
    \end{definition}

   \begin{definition} A Markov semigroup $\{P_t\}_{t\geq 0}$ on $B_b(\mathcal{X})$ is strongly continuous (see \cite{EK1986}) on $C_b(\mathcal{X})$, if for every $f\in C_b(\mathcal{X})$, 
			\begin{equation*}
				\lim\limits_{t\searrow 0}||P_tf-f||_{\infty}=0,
			\end{equation*}
   where ``$\searrow$'' or ``$\nearrow$'' stands for converging from above or below, respectively.
    \end{definition}

     \begin{definition}
        A Markov semigroup $\{P_t\}_{t\in T}$  is completely mixing (see \cite{J2013}), if  for every $f\in L_b(\mathcal{X})$ and for any $x,y\in\mathcal{X}$, 
	\begin{equation*}
	    \lim\limits_{t\rightarrow \infty}|P_tf(x)-P_tf(y)|=0.
	\end{equation*}
    \end{definition}
    
    \begin{definition}
        A Markov semigroup $\{P_t\}_{t\in T}$  is asymptotically stable, if  there exists a unique invariant measure $\mu_*\in\mathcal{P}(\mathcal{X}),$ and $P_t\mu $ converges weakly to $\mu_*$ for every $\mu\in\mathcal{P}(\mathcal{X})$ as $t\rightarrow\infty$.
    \end{definition}
    
    Comparing these notions, clearly, we may consider the following relations:
    \begin{itemize}
        \item  E-property $\underset{\mathbf{B_1}}{\overset{\mathbf{A_1}}{\rightleftharpoons}}$ Eventual e-property $\underset{\mathbf{B_2}}{\overset{\mathbf{A_2}}{\rightleftharpoons}}$ Eventual continuity;
        \item  Asymptotic stability $\underset{\mathbf{B_3}}{\overset{\mathbf{A_3}}{\rightleftharpoons}}$ Completely mixing property $\underset{\mathbf{B_4}}{\overset{\mathbf{A_4}}{\rightleftharpoons}}$ Eventual continuity.
    \end{itemize}\par 
    It can be checked that the implications $\mathbf{A_1}$ -  $\mathbf{A_4}$ follow form the definitions of these notions, which implies that the eventual continuity is a much weaker condition.    On the other hand,  generally, the implications $\mathbf{B_1}$ -  $\mathbf{B_4}$ are not satisfied. For example, 
    in \cite{J2013}, Jaroszewska provided sufficient conditions such that $\mathbf{B_4}$ holds for the eventually continuous Markov-Feller semigroups. In \cite{GL2015}, some criteria for
    the existence of invariant measures for the eventually continuous Markov-Feller semigroups are provided, and these results   also imply $\mathbf{B_3}$.  Moreover, in \cite{GLLL2023}, a necessary and sufficient condition for the asymptotic stability  is formulated directly for the eventually continuous Markov-Feller semigroups.\par

    In this paper, we mainly focus on $\mathbf{B_1}$ and $\mathbf{B_2}$. Given that both the eventual e-property and the eventual continuity do not imply the equicontinuity for the Markov-Feller semigroups, it is natural for us to wonder how to derive the e-property from these notions. We first show that the e-property and the eventual e-property are equivalent for discrete-time Markov-Feller semigroups (see Proposition \ref{Prop 1}). For the continuous-time case, we need to additionally assume the Markov-Feller semigroup is strongly continuous to ensure the equivalence (see Proposition \ref{Prop 2}).  Furthermore,  we provide a sufficient condition of the eventual e-property for the eventually continuous Markov-Feller semigroups. Thanks to the equivalence of the e-property and the eventual e-property, we hence obtain that the e-property and the eventual continuity are equivalent restricted on the support of each ergodic measure for the (discrete-time, or strongly continuous continuous-time) Markov-Feller semigroups (see Corollary \ref{Cor 1}).

	\section{Main results}\label{Sec 3}
	\subsection{Reation between the e-property and the eventual e-property}\label{Sec 3.1}
	We first show that the e-property and the eventual e-property are equivalent for the Markov-Feller semigroups.

	\begin{proposition}\label{Prop 1}
		Let $T=\mathbb{N}_+$. A Markov semigroup $\{P_t\}_{t\in T}$ satisfies the e-property if and only if it is a  Markov-Feller semigroup and satisfies the  eventual e-property. 
	\end{proposition}
	\begin{prof}[Proof.]
    ``$\Rightarrow$'': It suffices to show that the e-property implies the Feller property. By the definition of the e-property, it follows that $P_t(L_b(\mathcal{X})\subset C_b(\mathcal{X})$  which implies that  $P_t(C_b(\mathcal{X})\subset C_b(\mathcal{X})$ by \cite[Lemma 2.3]{HW2009} for all $t\in T$.\par 

    ``$\Leftarrow$'': It suffices to show that the e-property holds for all $x\in\mathcal{X}$. We prove it by contradiction.  Otherwise, assume that there exist $x\in \mathcal{X}$, $f\in L_b(\mathcal{X})$, $\epsilon>0$, $\{x_k\in\mathcal{X}:k\geq 1\}\rightarrow x$ and $\{n_k\in T:k\geq 1\}$ such that 

    \begin{equation*}
        \lim\limits_{k\rightarrow\infty} |P_{n_k}f(x_k)-P_{n_k}f(x)|=\epsilon>0.
    \end{equation*}
        Noting that $\{P_t\}_{t\in T}$ satisfies the eventual e-property, therefore $\{n_k:k\geq 1\}$ is  bounded. We may assume $1\leq n_k\leq N$ for some $N\in\mathbb{N}_+$ and all $k\geq 1.$ Due to the Feller property, it follows that 
	\begin{equation*}		0<\epsilon\leq\limsup\limits_{k\rightarrow\infty} |P_{n_k}f(x_k)-P_{n_k}f(x)|\leq \sum_{j=1}^{N}\limsup\limits_{k\rightarrow\infty} |P_{j}f(x_k)-P_{j}f(x)|=0,
		\end{equation*}
	which is a contradiction. We conclude that $\{P_t\}_{t\in T}$ satisfies the e-property.
	\end{prof}
	
	Moreover, this relation remains for the continuous-time Markov semigroups with the strong continuity on $C_b(\mathcal{X})$.

			\begin{proposition}\label{Prop 2}
			Let $T=\mathbb{R}_+$, and $\{P_t\}_{t\in T}$ be a Markov semigroup and strongly continuous on $C_b(\mathcal{X})$, then $\{P_t\}_{t\in T}$ satisfies the e-property if and only if it is Feller and satisfies the  eventual e-property. 
		\end{proposition}
		\begin{prof}[Proof.]
			We only need to show the opposite implication.	Assume that, contrary to our claim, the e-property fails at some $x\in\mathcal{X}$. Then there exist $f\in L_b(\mathcal{X})$, $\{x_k\in \mathcal{X}:k\geq 1\}\rightarrow x$, $\{t_k\in T:k\geq 1\}$  and $\epsilon>0$ such that 
			\begin{equation*}
			\lim\limits_{k\rightarrow\infty} |P_{t_k}f(x_k)-P_{t_k}f(x)|\geq \epsilon>0.
			\end{equation*}\par 
			Similar to Proposition \ref{Prop 1}, $\{t_k:k\geq 1\}$ is bounded by the eventual e-property. Passing by a subsequence if necessary, we may assume that $\{t_k:k\geq 1\}\searrow t_0\in\mathbb{R}_+$.	Let $g:=P_{t_0}f\in C_b(\mathcal{X})$ and $s_k=t_k-t_0$. Then we have 
			\begin{equation*}
			\lim\limits_{k\rightarrow\infty} |P_{s_k}g(x_k)-P_{s_k}g(x)|\geq \epsilon>0,
			\end{equation*}
			which conflicts the definition of the strong continuity on $C_b(\mathcal{X})$. We conclude that $\{P_t\}_{t\in T}$ satisfies the e-property.
		\end{prof}
		
    Indeed,  the eventual e-property not only deals with the ergodicity of the non-equicontinuous Markov-Feller semigroups (see \cite{C2017}), but also is helpful to handle the semigroup regularities of some SPDE models. We provide the next example to illustrate that the eventual e-property is a more convenient tool for SPDE models.

    \begin{example} \label{Ex heat equation}   Let $X$ satisfy the  stochastic heat equation on a torus   $\mathbb{T}=\mathbb{R}/\mathbb{Z}$:
		\begin{equation}\label{heat equation}
		d X(t,x)=\Delta X(t,x) dt + dW(t),\quad X(0,x)=\phi(x),\quad x\in\mathbb{T},\,t>0,
		\end{equation}
  where $W(t)$ will be determined later.
  Let $\{ e_k(x)=e^{{\rm i}kx}: x\in \mathbb{R},\;k\in \mathbb{Z}\}$ be a orthogonal basis of $L^2(\mathbb{T}):=\{\psi:\mathbb{T}\rightarrow \mathbb{R}:\int_{\mathbb{T}}\psi^2(x)dx<\infty\}$.  Let $W(t)$ be defined by $W(t)=\sum_{k\in\mathbb{Z}}\sigma_k e_k B_k(t)$, where $\sum_{k\in\mathbb{Z}} \sigma_k^2<\infty$, $\sigma_0=0$ and $\{B_k(t),\,t\in[0,\infty),\,k\in\mathbb{Z}\}$ are mutually independent real-valued standard Brownian motions on a probability space $(\Omega,\mathcal{F},\mathbb{P})$.  Let $L_0^2(\mathbb{T}):=\{ \psi\in L^2(\mathbb{T}): \int_{\mathbb{T}}\psi(x) dx =0\}$ and $ \phi\in L_0^2(\mathbb{T})$.  Moreover, let $H^1(\mathbb{T}):=\{ \psi\in L_0^2(\mathbb{T}): \int_{\mathbb{T}}(\frac{\partial \psi}{\partial x})^2dx <\infty\}$ 
 and $T(t):=e^{t\Delta}$ be the heat semigroup. By the same arguments as \cite[Theorem 5.4]{DPZ2014}, it follows that equation (\ref{heat equation}) has a unique weak solution in $L^2(\mathbb{T})$ given by the following formula
  \begin{equation}\label{heat equation2}
      X(t):=T(t)\phi+\int_{0}^{t}T(t-s)dW(s),\quad\text{for }t\geq0,\,\phi\in L^2(\mathbb{T}).
  \end{equation}  \par 
  Using integration by parts, we have that
  \begin{equation*}
      \langle X(t),e_k\rangle = \int_0^t\langle X(s),\Delta e_k\rangle ds + \sigma_k B_k(t),\quad \langle X(0), e_k\rangle = \langle \phi,e_k\rangle.
  \end{equation*}
    Let $X_k(t):=\langle X(t),e_k\rangle$ and $\phi_k:=\langle \phi,e_k\rangle$. Clearly, $X_k(t)$ satisfies the following O-U equation in one-dimension:
    \begin{equation*}
        X_k(t)=-k^2\int_0^tX_k(s)ds+\sigma_k B_k(t),\quad X_k(0) =  \phi_k,
    \end{equation*}
    hence, 
    \begin{equation*}
        X_k(t)= \phi_k e^{-k^2t}+\sigma_k\int_0^t e^{-k^2(t-s)}dB_k(s)\quad\text{for }k\in\mathbb{Z},
    \end{equation*}
    and
    \begin{equation*}
        X(t)=\sum_{k\in\mathbb{Z}}\phi_k e_k e^{-k^2t}+\sum_{k\in\mathbb{Z}}\sigma_k e_k\int_0^t e^{-k^2(t-s)}dB_k(s),\quad\text{for }t\geq 0.
    \end{equation*}\par 
  
  It can be checked that $X(t,\phi)\in L_0^2(\mathbb{H})$ almost surely for $t\geq 0$.  Let  $A(t)=\sum_{k\in\mathbb{Z}}\phi_k e_k e^{-k^2t}$ and $M(t)=\sum_{k\in\mathbb{Z}}\sigma_k e_k\int_0^t e^{-k^2(t-s)}dB_k(s)$. It follows that 
  \begin{equation}\label{heat equation3}
      ||A(t)||_{H^1}^2=\sum_{k\in\mathbb{Z}}\phi_k^2 k^2 e^{-2k^2t}\quad\text{and}\quad\mathbb{E}\int_0^t||M(t)||_{H^1}^2<\infty\quad\text{for }t\geq 0,
  \end{equation}
  which implies that $X(t,\phi)\in H^1(\mathbb{H})$ almost surely for any $t>0$ and $\phi\in L_0^2(\mathbb{H})$. In particular, equation (\ref{heat equation3}) implies that for any $\phi\in L_0^2(\mathbb{T})\setminus H^1(\mathbb{T})$, $||X(t,\phi||_{H^1}\nearrow\infty$ almost surely as $t\searrow 0$.\par 
  Moreover, it follows that 
 \begin{equation}\label{heat eqution4}
     ||X(t,\phi)-X(t,\tilde{\phi})||_{L^2}\leq e^{-t}||\phi-\tilde{\phi}||_{L^2}\text{ almost surely for } \phi,\tilde{\phi}\in L_0^2(\mathbb{T}),\, t\geq 0,
 \end{equation}
and that for any $s>0$,
\begin{equation}\label{heat eqution5}
     ||X(t,\phi)-X(t,\tilde{\phi})||_{H^1}\leq e^{-2(t-s)} ||X(s,\phi)-X(s,\tilde{\phi})||_{H^1} \text{ almost surely for } \phi,\tilde{\phi}\in L_0^2(\mathbb{T}),\, t\geq s.
 \end{equation}\par 
    Now we consider the Markov-Feller semigroup $\{P_t\}_{t\geq 0}$ generated by  equation (\ref{heat equation2}) on $\mathcal{X}:=L_0^2(\mathbb{T})$. More precisely, writing $X(t)=X(t,\phi)$, define
  \begin{equation*}
      P_tf(\phi)=\mathbb{E}f(X(t,\phi))\quad\text{ for }\phi\in \mathcal{X},\,f\in B_b(\mathcal{X}).
  \end{equation*}\par  
  Then the e-property in $\mathcal{X}$ follows from that for any $f\in L_b(\mathcal{X})$,
  \begin{equation*}
      |P_tf(\phi)-P_tf(\tilde{\phi})|\leq ||f||_{Lip} e^{-t}||\phi-\tilde{\phi}||_{L^2}\quad\text{for } \phi,\tilde{\phi}\in L_0^2(\mathbb{T}),\, t\geq 0.
  \end{equation*}\par 

  For the same reason, if we restrict $\{P_t\}_{t\geq 0}$ on $\mathcal{Y}:=H^1(\mathbb{T})$, i.e., 
    \begin{equation*}
      P_tf(\phi)=\mathbb{E}f(X(t,\phi))\quad\text{ for }\phi\in \mathcal{Y},\,f\in B_b(\mathcal{Y}),
  \end{equation*}
  then $\{P_t\}_{t\geq 0}$  still satisfies the e-property in $\mathcal{Y}$. However, if we let the starting point $\phi$ of equation (\ref{heat equation})  belong to $\mathcal{X}=L_0^2(\mathbb{T})$, then  $\{P_t\}_{t> 0}$ is still a Markov-Feller process on $\mathcal{Y}$.   In this case, we are unable to verify the e-property straightforwardly. For illustration, take $F(\varphi):=\sin (||\varphi||_{H^1}) \in L_b(\mathcal{Y})$. It follows that
\begin{equation*}
    \limsup\limits_{t\searrow 0} P_tF(\phi)=1\quad\text{and}\quad  \liminf\limits_{t\searrow 0} P_tF(\phi)=-1\quad\text{for any }\phi\in\mathcal{X}\setminus\mathcal{Y}.
\end{equation*}
Therefore,
\begin{equation*}
    \sup_{t>0}|P_tF(\phi)-P_tF(\tilde{\phi})|\geq 1/2\quad\text{for any }\phi\in\mathcal{X}\setminus\mathcal{Y},\,\tilde{\phi}\in\mathcal{Y},
\end{equation*}
which contradicts the e-property at $\tilde{\phi}$.\par  
  
   Instead, the eventual e-property in $\mathcal{Y}$ is satisfied by inequality (\ref{heat eqution5}). Therefore, we conclude that regardless of the starting points, the eventual e-property in $\mathcal{Y}$ is well-defined, but the e-property may cause trouble. Furthermore, the eventual e-property is useful to investigate the subtle long-time behaviors of Markov-Feller semigroups, for example,  where the topological support of the ergodic measure lives in.
   \end{example}

		\subsection{Equivalent condition for the e-property and the eventual continuity}\label{Sec 3.2}
	We are in a position to formulate the main result of our paper (Theorem \ref{Thm 1}). Indeed, we find out that for any Markov-Feller semigroup with an ergodic measure, and the interior of support of this ergodic measure is non-empty,  the e-property and the eventual continuity are equivalent on the interior of support.
        \begin{theorem}\label{Thm 1}
		Let  $\{P_t\}_{t\in T}$ be a  discrete-time, or strongly continuous continuous-time eventually continuous Markov-Feller semigroup, and let $\mu$ be an ergodic measure for $\{P_t\}_{t\in T}$.  If $\;\rm{Int}_{\mathcal{X}}(\rm{supp}\;\mu)\neq\emptyset$,  then $\{P_t\}_{t\in T}$ satisfies the e-property on $\rm{Int}_{\mathcal{X}}(\rm{supp}\;\mu)$.
		\end{theorem}

            Owing to Theorem \ref{Thm 1},  it directly shows that the eventual continuity and the e-property are equivalent, restricted on the support of each ergodic measure in its relative topology for any Markov-Feller semigroup.  For an ergodic measure $\mu$, denote $\mathcal{X}_{\mu}:=\text{supp}\;\mu$. As $\mathcal{X}_{\mu}$ is a closed set in $\mathcal{X}$, the subspace $(\mathcal{X}_{\mu},\rho)$ is still a  Polish space. 
			
	\begin{corollary}\label{Cor 1}  	Let  $\{P_t\}_{t\in T}$ be a  discrete-time, or strongly continuous continuous-time eventually continuous Markov-Feller semigroup, and let $\mu$ be an ergodic measure for $\{P_t\}_{t\in T}$.  Then  $\{P_t\}_{t\in T}$ has the e-property on $\mathcal{X}_{\mu}$.
		\end{corollary}
    \begin{prof}[Proof.]
        Given that $\mathcal{X}_{\mu}$ is an invariant set (see \cite[Lemma 4.1]{GL2015}), i.e.,
        \begin{equation*}
            P_t(x,\mathcal{X}_{\mu})=1\quad\text{for all } x\in\mathcal{X}_{\mu},\;t\in T,
        \end{equation*}
    hence $\{P_t\}_{t\in T}$ is still a Markov-Feller semigroup on $\mathcal{X}_{\mu}$. Clearly, $\rm{Int}_{\mathcal{X}_\mu}(\rm{supp}\;\mu)=\mathcal{X}_\mu$, which finishes the proof by Theorem \ref{Thm 1}.
        
    \end{prof}
       	
	In view of Proposition \ref{Prop 1}, \ref{Prop 2}, we reduce to verify the eventual e-property, which is much easier, to show the e-property for the Markov-Feller semigroups. In particular, we obtain a technical  criterion of the e-property, see Lemma \ref{Lemma 1}. The ideas follow from \cite{HSZ2017}.  

    \begin{lemma}\label{Lemma 1}
			  Let  $\{P_t\}_{t\in T}$ be a Markov-Feller semigroup and $\mu$ be an ergodic measure for $\{P_t\}_{t\in T}$. Assume that $\rm{Int}_{\mathcal{X}}(\rm{supp}\;\mu)\neq\emptyset$, then $\{P_t\}_{t\in T}$ satisfies the eventual e-property at $x_0\in\mathcal{X}$, if \\
			 (\textbf{a}) for any $f\in L_b(\mathcal{X})$, there exist a ball $B\subset \rm{supp}\; \mu$ and $N\in T$ such that
			\begin{equation}\label{eq Lemma1}
			|P_tf(x)-P_t f(x_0)|\leq\epsilon\;\;\text{for\;any\;}x\in B,\;t\geq N;
			\end{equation}
			 (\textbf{b}) There exist $\alpha>0$ and $t_0\in T$ such that $P_{t_0}\delta_{x_0}(B)>\alpha$ and ${\rm supp}\;P_{t_0}\delta_{x_0}\subset {\rm supp}\; \mu$. Moreover, for all $\nu\in\mathcal{P}({\rm supp}\;\mu),$ there exists $t_{\nu}\in T$ such that $P_{t_{\nu}}\nu(B)>\alpha$.\par 
    Condition (\textbf{b}) can also be replaced by:\\
                (\textbf{b'}) There exist $\alpha>0$ such that for all $\nu\in\mathcal{P}(\mathcal{X}),$ there exists $t_{\nu}\in T$ such that $P_{t_{\nu}}\nu(B)>\alpha$.
			\end{lemma}

    By Lemma \ref{Lemma 1} and Proposition \ref{Prop 1}, \ref{Prop 2}, it remains to check the conditions in Lemma \ref{Lemma 1} to obtain the e-property for the Markov-Feller semigroups. In order to prove Theorem \ref{Thm 1}, we need the following lemmas and the proofs are postponed in Section \ref{Sec 4}.
    
    \begin{lemma}\label{Lemma 4}
		Let $\mu,\nu\in\mathcal{P}(\mathcal{X})$ and $\mu$ be an invariant measure for $\{P_t\}_{t \in T}$. If $\rm{supp}\;\nu\subset \rm{supp}\;\mu$, then ${\rm supp}\;P_t\nu\subset \rm{supp}\;\mu$ for all $t\in T$.
	\end{lemma}
    
        \begin{lemma}\label{Lemma 3}
            Let  $\{P_t\}_{t\in T}$ a Markov-Feller semigroup. Assume that $\{P_t\}_{t\in T}$ is eventually continuous on $\mathcal{X}$ and $\mu$ is an ergodic measure for $\{P_t\}_{t\in T}$. Then for each $x\in\rm{supp}\;\mu,$ $Q_t(x,\cdot)$ weakly converges to $\mu$ as $t\rightarrow\infty$.
		\end{lemma}

         \begin{prof}[Proof of Theorem \ref{Thm 1}.]			
			Fix $x_0\in\rm{Int}_{\mathcal{X}}(\rm{supp}\;\mu)$, it suffices to show that the assumptions  in Lemma \ref{Lemma 1} hold at $x_0$. Condition (\textbf{a}) holds by the eventual continuity as follows. Fix $f\in L_b(\mathcal{X})$ and  $\epsilon>0.$  Due to the eventual continuity, there exists $\delta>0$ such that for $x\in B(x_0,\delta)\subset \rm{supp}\;\mu$,
			\begin{equation*}
				\limsup\limits_{t\rightarrow\infty}|P_tf(x)-P_tf(x_0))|\leq 2\epsilon.
			\end{equation*}
			Set $Y=\overline{B(x_0,\delta)}$  and $
			Y_{n}=\{x\in Y:|P_tf(x)-P_tf(x_0)|\leq\epsilon,\forall t\geq n\}$ for $n\in\mathbb{N}_+$.
			Note that $Y_{n}$ is closed and $Y=\bigcup_{n\geq 1} Y_{n}.$ By the Baire category theorem there exists $N\in\mathbb{N}$ such that Int$\;(Y_{N})\neq\emptyset$. Thus there exists some open set $B\subset Y_{N}$ such that 
			\begin{equation*}
			|P_tf(x)-P_t f(x_0)|\leq\epsilon\;\;\text{for\;any\;}x\in B,\;t\geq N. 
			\end{equation*}\par 
			Next we check Condition (\textbf{b}). Let $B(z,r)\subset B\subset \text{supp}\;\mu$. Then $\mu(B(z,r))>0$. Let $\alpha\in (0, \mu(B(z,r)))$, then by Lemma \ref{Lemma 3}, $Q_t(x_0,B)\rightarrow\mu(B)\geq \mu(B(z,r))>\alpha$ as $t\rightarrow\infty$. Hence there exists $t_0\in\mathbb{N}$ such that 
			$P_{t_0}\delta_{x_0}(B)>\alpha$.  Moreover, ${\rm supp}\;P_{t_0}\delta_{x_0}\subset {\rm supp}\;\mu$ by Lemma \ref{Lemma 4}. For any $\nu\in\mathcal{P}(\text{supp}\;\mu)$, 
			\begin{equation*}
			    \lim\limits_{t\rightarrow\infty}Q_t\nu(B)=\lim\limits_{t\rightarrow\infty}\int_{\text{supp} \;\mu}Q_t(x,B)\nu(dx)= \int_{\text{supp} \;\mu}\lim\limits_{t\rightarrow\infty}Q_t(x,B)\nu(dx)=\mu(B)>\alpha,
			\end{equation*}
			which implies that there exists $t_{\nu}\in T$, such that $P_{t_{\nu}}\nu(B)>\alpha$, completing the proof.
			\end{prof}

     As another application of Lemma \ref{Lemma 1}, we obtain \cite[Theorem 2.3]{HSZ2017} as a corollary.    We also notice that there is a small gap in \cite{HSZ2017},  where the definition of the e-property given in \cite{HSZ2017}  is  essentially equivalent to the eventual continuity (Definition \ref{def eventual continuity}). Therefore, the arguments in \cite[Theorem 2.3]{HSZ2017} indeed prove the eventual continuity for the asymptotically stable Markov-Feller semigroups. Actually, as illustrated in Section \ref{Sec 2}, the asymptotic stability itself ensures the eventual continuity without any other assumptions (the implications $\mathbf{A_3}$ and $\mathbf{A_4}$). We suppose that \cite[Theorem 2.3]{HSZ2017} is devoted to showing the e-property (Definition \ref{def e-property}), and the corresponding proofs can be revised by some modifications.

     \par

	\begin{corollary}
	  Let  $\{P_t\}_{t\in T}$ be an asymptotically stable (discrete-time, or strongly continuous continuous-time) Markov-Feller semigroup, and let $\mu_*$ be its unique invariant measure.    If $\;{\rm Int}_{\mathcal{X}}(\rm{supp}\;\mu_*)\neq\emptyset$, then $\{P_t\}_{t\in T}$ satisfies the e-property on $\mathcal{X}$.
	\end{corollary}

	\begin{prof}[Proof.]
		By Proposition \ref{Prop 1} and \ref{Prop 2}, it suffices to verify the assumptions in Lemma \ref{Lemma 1}. Fix $x_0\in\mathcal{X}$, Condition (\textbf{b'}) follows form the asymptotic stability, and  Condition (\textbf{a}) holds similar to \cite[Lemma 2.4]{HSZ2017} as follows.	Fix $f\in L_b(\mathcal{X}),\,\epsilon>0.$ Let $W$ be an open set such that $W\subset\text{\,supp}\;\mu_*.$ Set $Y=\overline{W}$ and 
			\begin{equation*}
			Y_{n}=\{x\in Y:|P_tf(x)-P_tf(x_0)|\leq\epsilon,\forall t\geq n\}\;\text{ for }n\in\mathbb{N}_+.
			\end{equation*}
			Note that $Y_{n}$ is closed and $Y=\bigcup_{n\geq 1} Y_{n}.$ By the Baire category theorem there exists $N\in\mathbb{N}$ such that Int$\;(Y_{N})\neq\emptyset$. Thus there exists some open set $B\subset Y_{N}$ such that 
			\begin{equation*}
			|P_tf(x)-P_tf(x_0)|\leq\epsilon\;\;\text{for\;any\;}x\in B,\;t\geq N.
			\end{equation*}
	\end{prof}

	\subsection{Discussion and conclusion}\label{Sec 3.3}

    Even though the results obtained in this paper characterize some relations between the e-property, the eventual e-property and the eventual continuity, a few questions are still left open.\par

    The first problem is how does the completeness of the space work for the semigroup regularities. For example, the strong continuity may fail if  the state space is not complete, but the e-property is still satisfied.  We have the next example for illustration.   Let $\mathcal{X}=[0,1],$ equipped with the Euclidean distance $\rho$. Let $S_t(x)= (x-t)_+=\max\{x-t,0\}$ and $P_tf(x)=f(S_t(x))=f((x-t)_+)$ for $t\geq 0$, $x\in\mathcal{X}$.   Clearly, $\{P_t\}_{t\geq 0}$ satisfies the e-property, the asymptotic stability and the strong continuity on $C_b(\mathcal{X})$. However, if we replace the metric $\rho$ by 
        \begin{equation*}
            d(x,y)=|x-y|,\quad\text{for }x,y\in[0,1),\quad\text{and}\quad d(x,1)=d(1,x)=1\quad\text{for } x\in[0,1),
        \end{equation*}
    then $\mathcal{X}$ is not complete. Moreover, the e-property and the asymptotic stability are still satisfied, but the strong continuity on $C_b(\mathcal{X})$ fails.\par

     The second problem is how to investigate the semigroup regularities for the Markov semigroups lacking the Feller property. We provide the following toy model to see that the e-property may fail for general Markov semigroups.  Let $\mathcal{X}=[-1,0]\cup\{1/n:n\geq 1\}$ equipped with the Euclidean distance $\rho$. Let $\gamma$ be a negative irrational number. Consider a deterministic dynamics $T:\mathcal{X}\rightarrow\mathcal{X}$ given by the following formula:
    \begin{equation*}
	T(x)=x+\gamma\mod 1\;\text{for}\;x\in[-1,0],\quad T(1/n)=1/(n-1)\;\text{for}\;n\geq 2, \quad \text{and}\quad T(1)=0.
    \end{equation*}	
    Let operator $P:\mathcal{M}(\mathcal{X})\rightarrow\mathcal{M}(\mathcal{X})$ be defined by $P\delta_x=\delta_{T(x)}$ for $x\in\mathcal{X}$. Clearly, the uniform distribution on $[-1,0]$ is the unique invariant measure. Moreover, the eventual continuity holds on $\mathcal{X}$, and the e-property holds on $\mathcal{X}\setminus\{0\}$. Yet the (eventual) e-property  fails at $0$, and this  Markov semigroup does not satisfy the Feller property. 
    Due to Proposition \ref{Prop 1}, the Feller property is necessary for the e-property, so the e-property fails for the Markov semigroups lacking the Feller property. On the other hand, however, the eventual continuity may hold for  general Markov semigroups lacking the Feller property, similar to this example. Therefore, it seems likely that the eventual continuity is a suitable tool to deal with the ergodicity for the Markov semigroups lacking the Feller property.\par

 Thirdly, although Theorem \ref{Thm 1} provides a general approach to verify the e-property, it requires the Markov-Feller semigroup to possess the eventual continuity at each point on the whole space. Practically, this assumption may be too demanding in some circumstances, when we are only able to verify the eventual continuity for only one point or several points instead of the whole space. In these cases, how can we obtain the e-property?    For example, \cite[Theorem 1]{GLLL2023} shows that, provided that  both the eventual continuity and a lower bound condition hold at one certain point, the  asymptotic stability follows, which can be further applied to  guarantee the e-property by Corollary \ref{Cor 1}.   We conjecture that if the eventual continuity is satisfied for one point  in the support of an ergodic measure, then the  eventual continuity also holds on the whole support, and hence implies the e-property (on the support). \par

    The next problem is how to show the \emph{Ces\`aro e-property}, which is defined for the Ces\`aro averages of the Markov-Feller semigroups. This notion was formulated by Worm in  \cite{W2010}, which is a generalization of the e-property. Worm proved several applications on the Yosida-type ergodic decomposition of state space for the Ces\`aro e-property in \cite{W2010}. The Ces\`aro e-property is weaker than the e-property and is not generally satisfied. For example, Hile et al. in \cite{HSZ2017} provided an asymptotically stable Markov-Feller semigroup which does not satisfy the Ces\`aro e-property. One natural thought is to apply the \emph{Ces\`aro eventual continuity} (see definition in \cite{GL2015}), which generalizes the eventual continuity, to show the Ces\`aro e-property as in Theorem \ref{Thm 1}. However, as the techniques used in Theorem \ref{Thm 1} are inapplicable for the Ces\`aro averages, there is still no proper way to  show the Ces\`aro e-property by means of the Ces\`aro eventual continuity as far as we know.\par

    We also expect to figure out the relations between the time-continuity of the continuous-time Markov-Feller semigroups and the semigroup regularities.  In Proposition \ref{Prop 2}, we assume that the Markov-Feller semigroups are strongly continuous on $C_b(\mathcal{X})$, which is a rather strong assumption and is not generally satisfied. Therefore, we wonder whether we can use some weaker assumptions on the time-continuity to maintain the equivalence between the e-property and the eventual e-property for the Markov-Feller semigroups.  For example, whether can we replace the strong continuity on $C_b(\mathcal{X})$,  by the stochastic continuity on $C_b(\mathcal{X})$, i.e.,
    $P_tf(x)\rightarrow f(x)$ as $t\searrow 0$ for all $x\in\mathcal{X}$, $f\in C_b(\mathcal{X})$?  Or whether can we replace the strong continuity on $C_b(\mathcal{X})$ by a smaller space, $C_0(\mathcal{X})$, i.e. the space of continuous functions which vanish at infinity? For illustration,  we consider the state space $\mathcal{X}$ to be a locally compact separable metric space for now, which is still a Polish space. The strong continuity on $C_0(\mathcal{X})$ is equivalent to the stochastic continuity on $C_0(\mathcal{X})$ (see \cite{EK1986}).  In this case, it suffices to guarantee stochastic stability on $C_0(\mathcal{X})$ to obtain equivalence between the e-property and the eventual e-property in $C_b(\mathcal{X})$, once we can 
    replace the strong continuity assumption on $C_b(\mathcal{X})$ by $C_0(\mathcal{X})$. \par

    Another interesting problem is the relation between the \emph{asymptotic strong Feller} property and the e-property (or the eventual continuity). The asymptotic strong Feller property is a generalization of the strong Feller property, and was formulated by Hairer and Mattingly to attack the unique ergodicity for 2D Navier-Stokes equations with highly degenerate noise in \cite{HM2006,HM2008,HM2011}. The definition of the asymptotic strong Feller property is rather involved, and we refer the readers to \cite{HM2006}. Generally, the asymptotic strong Feller property and the e-property do not imply each other. Furthermore, the asymptotic strong Feller property does not imply the eventual continuity and the strong continuity. Noting that the lack of time-regularity may break the semigroup regularities, we wonder whether  we can deduce better semigroup regularities by promoting the time-continuity.  We conjecture that the asymptotic strong Feller property together with the strong continuity imply the eventual continuity.   For illustration, we provide the next example which neither satisfies the strong continuity on $C_b(\mathcal{X})$, nor the eventual continuity (and the e-property), but it satisfies the asymptotic strong Feller property.

    \begin{example}\label{example}
            
    The construction is based on the Jaroszewska's ideas in \cite{J20132}. A set $\mathbb{B}\subset\mathbb{R}$ is called a Hamel basis for $\mathbb{R}$ if every element of $\mathbb{R}$ is a unique finite rational linear combination of elements of $\mathbb{B}$. The existence of a Hamel basis is guaranteed by the axiom of choice. The Hamel bases are useful for constructions of functions with nontypical properties, as the following lemma, taken from \cite{RS1998}.
		\begin{lemma}
            [\cite{RS1998}, Theorem 1.6]\label{Lemma ex}
			If $\mathbb{B}$  is a Hamel basis for $\mathbb{R}$ and $g:\mathbb{B}\rightarrow \mathbb{R}$ is an arbitrary function then there exists a function
			$\varphi_g:\mathbb{R}\rightarrow \mathbb{R}$ which satisfies the Cauchy equation (i.e., $\varphi_g(x)+\varphi_g(y)=\varphi_g(x+y)$ for all $x,y\in\mathbb{R}$, in other words, $\varphi_g$ is additive) and such that $\varphi_g|_{\mathbb{B}}=g|_{\mathbb{B}}$.
		\end{lemma}
  
		Let $\mathcal{X}=\mathbb{R}$. Fix a sequence $\{b_n\in\mathbb{B}:n\geq 1\}$ and let $b_0:=0$. Then $b_i\neq b_j$ for $i\neq j$ by definition. Moreover, there exists $\alpha_n\in\mathbb{Q}$ for $n\in\mathbb{N}$ such that $\alpha_{n}b_{n}+1\leq\alpha_{n+1} b_{n+1}\leq \alpha_{n}b_{n}+2$. Let $\widetilde{\mathbb{B}}:=(\mathbb{B}\setminus\{b_n:n\geq 1\})\cup \{\alpha_n b_n:n\geq 1\}$, and $\widetilde{\mathbb{B}}$ is also a Hamel basis for $\mathbb{R}$. Define $a_n:=\alpha_n b_n$ for $n\in\mathbb{N}_+$. Then $0<a_n\leq a_1+2n$ for $n\in\mathbb{N}_+$ and $a_n\nearrow\infty$ as $n\rightarrow\infty$.	 Let $g:\widetilde{\mathbb{B}}\rightarrow\mathbb{R}$ such that $g(a_{2k-1})=k^3$ and $g(a_{2k})=-k^3$ for $k\geq 1$. Let $\varphi_g$ be the additive extension of $g$ by Lemma \ref{Lemma ex}. Next, let
		\begin{equation*}
			S_t(x)=e^{\varphi_g(a_1t)}x\quad\text{ for }t\geq 0,x\in\mathcal{X}.
		\end{equation*}
        As $\varphi_g$ is additive, $\{S_t\}_{t\geq 0}$ is a semigroup and hence it generates a semigroup  $\{P_t\}_{t\geq 0}$ by $P_t\delta_x:=\delta_{S_t(x)}$ for $t\geq 0$, $x\in\mathcal{X}$. Since $S_t$ is continuous for each $t\geq 0$, $\{P_t\}_{t\geq 0}$ is a Markov-Feller semigroup, with the dual semigroup, still denote by $\{P_t\}_{t\geq 0}$, given by
        \begin{equation*}
            P_tf(x)=f(S_t(x))=f(e^{\varphi_g(a_1t)}x),\quad\text{for } t\geq 0,\,x\in\mathcal{X},\,f\in B_b(\mathcal{X}).
        \end{equation*} 
	   Take  $t_k=k^{-2}a_1^{-1}a_{2k-1}$ and $s_k=k^{-2}a_1^{-1}a_{2k}$ for $k\geq 1$. By definition, $t_k,\;s_k\searrow 0$ as $k\rightarrow\infty$. Fix $x>0$, by the additivity of $\varphi_g$,
		\begin{equation*}
			S_{t_k}(x)=e^{\varphi_g(a_1t_k)}x=e^{k^{-2}\varphi_g(a_{2k-1})}x=e^{k^{-2}g(a_{2k-1})}x=e^kx\rightarrow\infty\quad\text{as }k\rightarrow\infty.
		\end{equation*}
		For the same reason, $S_{s_k}(x)\rightarrow-\infty$ as $k\rightarrow\infty$, which in turn, conflicts the strong continuity on $C_b(\mathcal{X})$.\par 
        Next we examine the eventual continuity at $0\in\mathcal{X}$. Let $f(\cdot)=|\cdot| \wedge 1\in L_b(\mathcal{X})$. Then for any $ y\in B(0,1/2)\setminus\{0\}$, 
        \begin{equation*}
            \limsup\limits_{t\rightarrow\infty} |P_tf(y)-P_tf(0)|\geq \limsup\limits_{t\rightarrow\infty,t\in\{\frac{a_{2k-1}}{a_1}:k\in\mathbb{N}_+\}} |P_tf(y)-P_tf(0)|= \limsup\limits_{k\rightarrow\infty} (e^{k^3}|y| \wedge 1)=1,
        \end{equation*}
    hence the eventual continuity fails at $0\in\mathcal{X}$.\par 
    Finally we show that $\{P_t\}_{t\geq 0}$ satisfies the asymptotic strong Feller property for any $x\in\mathcal{X}$ with $t_k=a_1^{-1}a_{2k}$ and $\rho_k= 1\wedge k\rho$. Then the asymptotic strong Feller property at $x\in\mathcal{X}$ follows from that
    \begin{equation*}
    \begin{aligned}
        ||P_{t_k}\delta_x-P_{t_k}\delta_y||_{\rho_k}&\leq \sup_{f\in L_b(\mathcal{X}), ||f||_{Lip}\leq 1} k|P_{t_k}f(x)-P_{t_k}f(y)|\\
        &\leq k e^{\varphi_g(a_1t_k)}\rho(x,y)=ke^{-k^3}|x-y|\rightarrow 0\quad\text{as }k\rightarrow\infty,
      \end{aligned}
    \end{equation*}
    where $||\mu-\nu||_{d}$ is the Wasserstein distance between $\mu$, $\nu\in\mathcal{P}(\mathcal{X})$ with respect to the metric $d$ on $\mathcal{X}$. 
    \end{example}

	\section{Proofs of Lemmas}\label{Sec 4}

	\begin{prof}[Proof of Lemma \ref{Lemma 1}.]	
         Assume that, contrary to our claim, $\{P_t\}_{t\in T}$ does not satisfy the eventual e–property at $x_0$. Therefore there exists a function $f\in L_b(\mathcal{X})$  such that
		\begin{equation*}
		\limsup\limits_{(x,t)\rightarrow(x_0,\infty)}|P_tf(x)-P_tf(x_0)|>0.
		\end{equation*}
		We may choose $\epsilon>0$ and $x_j\rightarrow x_0,t_j\rightarrow\infty$ as $j\rightarrow\infty$ such that
		\begin{equation*}
		\limsup\limits_{j\rightarrow\infty}|P_{t_j}f(x_j)-P_{t_j}f(x_0)|>3\epsilon.
		\end{equation*}		
		Let $B:=B(z,r)$ be a ball such that condition (\ref{eq Lemma1}) holds. 
		Let $k\geq 1$ be such that $2(1-\alpha)^k|f|_{\infty}<\epsilon.$ By induction we are going to define two sequences of measures $\{\nu_i^{x_0}\}_{i=1}^k,\{\mu_i^{x_0}\}_{i=1}^k,$ and a sequence of positive numbers $\{s_{i}\}_{i=1}^k$ in the following way: by Condition (\textbf{b}) or (\textbf{b'}), let
		$s_1>0$ be such that	
		\begin{equation*}
		P_{s_1}\delta_{x_0}(B(z,r))>\alpha.
		\end{equation*}\par
		Choose $r_1<r$ such that $P_{s_1}\delta_{x_0}(B(z,r_1))>\alpha$ and $P_{s_1}\delta_{x_0}(\partial B(z,r_1))=0$ and set
		\begin{center}
			$\nu_1^{x_0}(\cdot) = \dfrac{P_{s_1}\delta_{x_0}(\cdot \cap B(z,r_1))}{P_{s_1}\delta_{x_0}(B(z,r_1))},\quad$
			$\mu_1^{x_0}(\cdot) = \dfrac{1}{1-\alpha}(P_{s_1}\delta_{x_0}(\cdot)-\alpha\nu_1^{x_0}(\cdot)).$
		\end{center}
        If Condition (\textbf{b}) holds, then $P_{s_1}\delta_{x_0}\subset \text{supp}\;\mu$ and $\text{supp}\;\nu_1^{x_0}\subset B(z,r)\subset \text{supp}\;\mu$, hence $\text{supp}\;\mu_1^{x_0}\subset \text{supp}\;\mu$.\par 
		Assume that we have done it for $i = 1,\dots , l,$ for some $l < k.$ By Condition (\textbf{b}) or (\textbf{b'}), now let $s_{l+1}$ be such that
		$P_{s_{l+1}}\mu_l^{x_0}(B(z,r))>\alpha$. 	Choose $r_{l+1}<r$ such that $P_{s_{l+1}}\mu_l^{x_0}(B(z,r_{l+1}))>\alpha$ and $P_{s_{l+1}}\mu_l^{x_0}(\partial B(z,r_{l+1}))=0$ and set
		
		\begin{center}
			$\nu_{l+1}^{x_0}(\cdot) = \dfrac{P_{s_{l+1}}\mu_l^{x_0}(\cdot \cap B(z,r_{l+1}))}{	P_{s_{l+1}}\mu_l^{x_0}(B(z,r_{l+1}))},\quad$
			$\mu_{l+1}^{x_0}(\cdot) = \dfrac{1}{1-\alpha}(P_{s_{l+1}}\mu_l^{x_0}(\cdot)-\alpha\nu_{l+1}^{x_0}(\cdot)).$
		\end{center}
		Then it gives 
		\begin{equation*}
		\begin{aligned}
		P_{s_1+\dots+s_k}\delta_{x_0}(\cdot)&=\alpha P_{s_2+\dots+s_k}\nu_1^{x_0}(\cdot)+\alpha(1-\alpha) P_{s_3+\dots+s_k}\nu_2^{x_0}(\cdot)+ \cdots \\
		&+\alpha(1-\alpha)^{k-1}\nu_k^{x_0}(\cdot)+(1-\alpha)^k \mu_k^{x_0}(\cdot),
		\end{aligned}
		\end{equation*}
		similarly, $\text{supp}\;\nu_{l+1}^{x_0}\subset B(z,r)\subset \text{supp}\;\mu_*$ and $\text{supp}\;\mu_{l+1}^{x_0}\subset \text{supp}\;\mu_*$. \par 
		We further adopt the same procedure to construct the sequence $\{\nu_i^{x_j}\}_{i=1}^k,\{\mu_i^{x_j}\}_{i=1}^k$ for each $j\geq 1$.
		And by the same arguments as in \cite{HSZ2017}, it turns out that $\nu_i^{x_j}$ converges weakly to $\nu_i^{x_0}$ and 
		$\mu_i^{x_j}$ converges weakly to $\mu_i^{x_0}$  as $j\rightarrow\infty$ for $i=1,\dots,k$.\par 
		Observe that for any $x_j$ sufficiently close to $x_0$ and all $t\geq s_1+\cdots+s_k, $ we have
		\begin{equation*}
		\begin{aligned}
		P_t\delta_{x_j}(\cdot)= &\,\alpha P_{t-s_1}\nu_1^{x_j}(\cdot)+\alpha(1-\alpha) P_{t-s_1-s_2}\nu_2^{x_j}(\cdot)+ \cdots \\
		&+\alpha(1-\alpha)^{k-1}P_{t-s_1-\cdots-s_k}\nu_k^{x_j}(\cdot)+(1-\alpha)^k P_{t-s_1-\cdots-s_k}\mu_k^{x_j}(\cdot),
		\end{aligned}
		\end{equation*}
		where  $\text{supp }{{\nu}_i^{x_j}}\subset B(z,r),\,j\geq 1,\,i=1,\dots,k.$\par 
		Therefore, by (\ref{eq Lemma1}), 
		\begin{equation*}
		|\langle P_tf-P_tf(x_0),\nu_i^{x_j}\rangle|\leq\int_{\mathcal{X}} |P_tf(y)-P_tf(x_0)|\nu_i^{x_j}(dy)\leq \epsilon/2\quad\text{for all}\;j\geq 1,\,i=1,\dots,k,\,t\geq N.
		\end{equation*}
		The same inequality also holds for $\nu_i^{x_0},i=1,\dots,k$. Thus it follows that
		\begin{equation*}
		|\langle f,P_t\nu_i^{x_j}\rangle-\langle f,P_t\nu_i^{x_0}\rangle|=|\langle P_tf-P_tf(x_0),\nu_i^{x_j}\rangle-\langle P_tf-P_tf(x_0),\nu_i^{x_0}\rangle|\leq\epsilon,
		\end{equation*}
	for all $j\geq 1$, $i=1,\dots,k$, $t\geq N$. Furthermore, we obtain that
		\begin{equation*}
		\limsup\limits_{j\rightarrow\infty}|\langle f,P_{t_j}\nu_i^{x_j}\rangle-\langle f,P_{t_j}\nu_i^{x_0}\rangle|\leq\epsilon\quad\text{for}\;i=1,\dots,k.
		\end{equation*}
		Hence, it follows that
		\begin{equation*}
		\begin{aligned}
		3\epsilon&<\limsup\limits_{j\rightarrow\infty}|P_{t_j}f(x_j)-P_{t_j}f(x_0)|=\limsup\limits_{j \to \infty}|\langle f,P_{t_j} \delta_{x_j}\rangle-\langle f,P_{t_j} \delta_{x_0}\rangle| \\
		&\leq \alpha\limsup\limits_{j \to \infty}|\langle f, P_{t_j}\nu_1^{x_j} \rangle-\langle f, P_{t_j}\nu_1^{x_0} \rangle|+\cdots+\alpha(1-\alpha)^l\limsup\limits_{j \to \infty}|\langle f, P_{t_j}\nu_l^{x_j}\rangle-\langle f, P_{t_j}\nu_l^{x_0} \rangle|\\
        &\quad+2(1-\alpha)^k|f|_{\infty}\\
		&\leq (\alpha+\cdots+\alpha(1-\alpha)^l)\epsilon+\epsilon\\
		&\leq 2\epsilon,
		\end{aligned}
		\end{equation*}
		which is impossible. This completes the proof.
	\end{prof}

	\begin{prof}[Proof of Lemma \ref{Lemma 4}.]
		It suffices to prove for $\nu=\delta_x$, $x\in\rm{supp}\;\mu$. Fix $t\in T$  and $\epsilon>0$. Let $z\in {\rm supp}\;P_t\delta_x$. Then $P_t\delta_{x}(B(z,\epsilon))>0$. Moreover, due to the Feller property, there exists $\delta>0$ such that $P_t\delta_{y}(B(z,\epsilon))> \frac{1}{2}P_t\delta_{x}(B(z,\epsilon))>0$ for all $y\in B(x,\delta)$.
		Since $x\in \rm{supp}\;\mu$, it follows that
		\begin{equation*}
		\begin{aligned}
		\mu(B(z,\epsilon))&=\int_{\mathcal{X}}P_t(y,B(z,\epsilon))\mu(dy)\geq \int_{B(x,\delta)}P_t(y,B(z,\epsilon))\mu(dy)\\
    &\geq \mu(B(x,\delta)) \cdot \inf_{y\in B(x,\delta)} P_t(y,B(z,\epsilon))>0,
		\end{aligned}
		\end{equation*}
		which completes the proof.
	\end{prof}\par

	\begin{prof}[Proof of Lemma \ref{Lemma 3}.] The proof is divided into three parts. \\
	$\mathbf{Step\;1.}$ We first show that 	$\{Q_t(x,\cdot)\}_{t\in T}$ is tight for any $x\in\text{supp\;}\mu$. 
	Assume that, on the contrary, $\{Q_t(x,\cdot)\}_{t\in T}$ is not tight for some $x\in\text{supp\;}\mu$. Then by \cite[Lemma 1]{KPS2010}, there exist a strictly increasing sequence of positive numbers $t_i\nearrow\infty$, a positive number $\epsilon$ and a sequence of compact sets $\{K_i\}$ such that
		\begin{equation}\label{eq 5.5}
		Q_{t_i}(x,K_i)\geq\epsilon,\quad\forall i,\quad\text{and}\quad \min\{\rho(x,y):x\in K_i,y\in K_j\}\geq\epsilon,\quad\forall i\neq j.
		\end{equation}
		We will derive the assertion from the claim that there exist sequences $\{\bar{f}_k\}\subset L_b(\mathcal{X})$, $\{\nu_k\}\subset\mathcal{P}(\mathcal{X})$ and an increasing sequence of integers $\{m_k\}$ such that supp $\nu_k\subset B(x,1/k)$ for $k\in\mathbb{N}_+$, and
		\begin{equation}\label{eq 5.7}
		\mathbf{1}_{K_{m_k}}\leq\bar{f}_k\leq\mathbf{1}_{K_{m_k}^{\epsilon/4}}\quad\text{and}\quad \text{Lip}(\bar{f}_k)\leq 4/\epsilon;
		\end{equation}
		\begin{equation}\label{eq 5.8}
		Q_t\nu_k(\cup_{i=k}^{\infty}K_{m_i}^{\epsilon/4})\leq\epsilon/4,\quad\text{for all } t\in T;
		\end{equation}
		\begin{equation}\label{eq 5.9}
		\limsup\limits_{t \to \infty}|\langle Q_t\delta_{x},f_k\rangle-\langle Q_t\nu_k,f_k\rangle|\leq\epsilon/4,
		\end{equation}
		where $f_1:=0,f_k:=\sum_{i=1}^{k-1}\bar{f}_i,n\geq 2.$  Let $f:=\sum_{i=1}^{\infty}\bar{f}_i$. By (\ref{eq 5.5}) and (\ref{eq 5.7}), $f$ is uniformly bounded with $||f||_{\infty}=1.$ Further, noting that for any $x,y\in X$ with $\rho(x,y)<\epsilon/8,$ we have $\bar{f}_i(x)\neq0,$ or $\bar{f}_i(y)\neq0$ for at most one $i.$  Thus
		\begin{equation*}
		|f(x)-f(y)|\leq 16\epsilon^{-1} \rho(x,y)
		\end{equation*} 
		and $f\in L_b(\mathcal{X}).$ Then it follows that
		\begin{equation}\label{eq 5.10}
		\begin{aligned}
		\langle Q_t\delta_{x},f\rangle -\langle Q_t\nu_k,f\rangle\geq Q_t(x,\cup_{i=k}^{\infty}K_{m_i})+\langle Q_t\delta_{x},f_k\rangle
		-\langle Q_t\nu_k,f_k\rangle-Q_t\nu_k(\cup_{i=k}^{\infty}K_{m_i}^{\epsilon/4}).
		\end{aligned}
		\end{equation}
		By (\ref{eq 5.5}),
		\begin{equation}\label{eq 5.11}
		\limsup\limits_{t \to \infty}Q_t(x,\cup_{i=k}^{\infty}K_{m_i})\geq	\limsup\limits_{l \to \infty}Q_{t_{m_l}}(x,\cup_{i=k}^{\infty}K_{m_i})\geq\limsup\limits_{l \to \infty}Q_{t_{m_l}}(x,K_{m_l})\geq\epsilon.
		\end{equation}
		From (\ref{eq 5.8})-(\ref{eq 5.11}),  it follows that
		\begin{equation*}
		\limsup\limits_{t \to \infty}[\langle Q_t\delta_{x},f\rangle -\langle Q_t\nu_k,f\rangle]\geq\epsilon-\epsilon/4-\epsilon/4=\epsilon/2.
		\end{equation*}
		Hence there must be a sequence $y_k\in\text{supp\;}\nu_k$ such that
		\begin{equation*}
		\limsup\limits_{t \to \infty}|Q_tf(x)-Q_tf(y_k)|\geq\epsilon/2,
		\end{equation*}
		which contradicts the eventual continuity of $\{P_t\}_{t \in T}$ at $x.$ This completes the proof.\par 
		$Proof\;of\;the\;claim.$ We accomplish this by induction on $k$. Let $k = 1$. Given $x\in\text{supp\;}\mu$, we have $\mu(B(x,\delta))>0$ for all $\delta>0$. Let $\nu_1\in\mathcal{P}(\mathcal{X})$ be defined by the formula
		\begin{equation*}
		\nu_1(B)=\mu(B|B(x,1)):=\frac{\mu(B\cap B(x,1))}{\mu(B(x,1))},\quad B\in\mathcal{B}(\mathcal{X}).
		\end{equation*}
		Since $\nu_1\leq\mu^{-1}(B(x,1))\mu$, from the fact that $\mu$ is ergodic, it follows that the family $\{Q_t\nu_1\}_{t\in T}$ is tight. Then there exists some compact set $K$ such that 
		\begin{equation*}
		Q_t\nu_1(K^c)\leq\epsilon/4,\quad\text{for all }t\in T.
		\end{equation*} 
		Note, however, that $K\cap K_{i}^{\epsilon/4}\neq\emptyset$ for only finitely many $i'$s. As a result, there exists an integer $m_1$ such that
		\begin{equation*}
		Q_t\nu_1(\cup_{i=1}^{\infty}K_{m_1}^{\epsilon/4})\leq\epsilon/4,\quad\text{for all } t\in T.
		\end{equation*}
		Let $\bar{f}_1$ be an arbitrary Lipschitz function satisfying
		\begin{equation*}
		\mathbf{1}_{K_{m_1}}\leq\bar{f}_1\leq\mathbf{1}_{K_{m_1}^{\epsilon/4}}\quad\text{and}\quad \text{Lip}(\bar{f}_1)\leq 4/\epsilon.
		\end{equation*}
		Assume, now, that for a given $k\geq 1$, we have already constructed $\bar{f}_1,\dots,\bar{f}_k$, $\nu_1,\dots,\nu_k$ and $m_1,\dots,m_k$ satisfying the claim. In view of the  eventual continuity of $\{P_t\}_{t \in T}$, we can choose $\delta\in(0,1/(k+1))$ such that 
		\begin{equation*}
		\sup\limits_{y\in B(x,\delta)}\limsup\limits_{t \to \infty}|Q_tf_{k+1}(x)-Q_tf_{k+1}(y)|<\epsilon/4.
		\end{equation*}
		Further, let $\nu_{k+1}(\cdot):=\mu(\cdot|B(x,\delta)).$ Therefore, by the dominate convergence theorem,
		\begin{equation*}
		\limsup\limits_{t \to \infty}|\langle Q_t\delta_{x},f_{k+1}\rangle-\langle Q_t\nu_{k+1},f_{k+1}\rangle|\leq\epsilon/4.
		\end{equation*}
		Finally, we let $\bar{f}_{k+1}$ be an arbitrary bounded, globally Lipschitz function satisfying (\ref{eq 5.7}).\par

	$\mathbf{Step\;2.}$ Next, we show that 	$\{Q_t(x,\cdot)\}_{t\in T}$ weakly converges to some invariant measure $\nu$ for each $x\in\text{supp }\mu$. 	Assume that, on the contrary, there exists  $x\in\text{supp }\mu$ such that the sequence 	$\{Q_t(x,\cdot)\}_{t\in T}$ does not converge. Given that 	$\{Q_t(x,\cdot)\}_{t\in T}$ is tight, by the Prokhorov theorem we may find at least two different probability measures $\mu_1,\mu_2$ and two sequences   $\{s_n\}_{n\geq 1}\nearrow\infty,\;\{t_n\}_{n\geq 1}\nearrow\infty$ such that $\{Q_{s_n}(x,\cdot)\}_{n\geq 1},\;\{Q_{t_n}(x,\cdot)\}_{n\geq 1}$ weakly converges to $\mu_1,\mu_2,$ respectively.\par 
		Choose $f\in L_b(\mathcal{X})$ and $\epsilon>0$ such that $|\langle\mu_1,f\rangle-\langle\mu_2,f\rangle|>\epsilon.$ \par 
		Let $D=\{x_k\in \mathcal{X}:k\geq 1\}$ be a countable dense set of $\mathcal{X}$. Passing to a subsequence if necessary, we may assume that $\lim\limits_{n\rightarrow\infty}Q_{s_n}f(x_k)$ exists for all $k\geq 1.$ Now let $\overline{g}(x):=\limsup\limits_{n\rightarrow\infty}Q_{s_n}f(x)$ and $\underline{g}(x):=\liminf\limits_{n\rightarrow\infty}Q_{s_n}f(x)$ for $x\in \mathcal{X}.$ We claim that $\overline{g}=\underline{g}\in C_b(\mathcal{X}),$ and denote $g=\overline{g}=\underline{g}.$ Actually, by the eventual continuity of $\{P_t\}_{t\in T}$, for $x\in \mathcal{X},\;\eta>0,$ there exists some $x_k$ such that
		\begin{equation*}
		\limsup\limits_{n\rightarrow\infty}|Q_{s_n}f(x)-Q_{s_n}f(x_k)|\leq\eta/2,
		\end{equation*}
		we have
		\begin{equation*}
		\overline{g}(x)-\underline{g}(x_k)\leq\eta/2\quad\text{and }\quad\underline{g}(x)-\overline{g}(x_k)\leq\eta/2,
		\end{equation*}
		hence
		\begin{equation*}
		|\overline{g}(x)-\underline{g}(x)|\leq\eta,\quad\text{for all }x\in \mathcal{X},\eta>0.
		\end{equation*}
		Thus we conclude that $\{Q_{s_n}f\}_{n\geq 1}$ converges to $g\in C_b(\mathcal{X})$ pointwisely. \par 
		
		By the bounded convergence theorem and  invariance, we have
		\begin{equation*}
		\langle \mu_2,f\rangle =\lim\limits_{n\rightarrow\infty}\langle \mu_2, Q_{s_n}f\rangle=\langle\mu_2,g\rangle.
		\end{equation*}
		As $\{Q_{t_n}(x,\cdot)\}_{n\geq 1}$ converges weakly to $\mu_{2},$ we can fix $N\in\mathbb{N}_+$ such that
		\begin{equation}\label{eq 5.1}
		|\langle\mu_{2},g\rangle-\langle Q_{t_N}(x,\cdot),g\rangle|\leq \epsilon/5.
		\end{equation}
		For such $N,$ we  choose $n$ sufficiently large such that
		\begin{equation}\label{eq 5.2}
		|\langle Q_{t_N}(x,\cdot),g\rangle-\langle Q_{t_N}(x,\cdot),Q_{s_n}f\rangle|\leq \epsilon/5
		\end{equation}
		Further, by \cite[Lemma 2]{KPS2010}, we have
		\begin{equation*}
		\lim\limits_{n\rightarrow\infty}||Q_{s_n,t_N}(x,\cdot)-Q_{s_n}(x,\cdot)||_{TV}=0.
		\end{equation*}
		Hence we fix $n$ sufficiently large such that
		\begin{equation}\label{eq 5.3}
		|\langle Q_{s_n,t_N}(x,\cdot),f\rangle-\langle Q_{s_n}(x,\cdot),f\rangle|\leq \epsilon/5
		\end{equation}
		Finnaly, note that $\{Q_{s_n}(x,\cdot)\}_{n\geq 1}$ converges weakly to $\mu_{1},$ there exists $n$ sufficiently large such that
		\begin{equation}\label{eq 5.4}
		|\langle Q_{s_n}(x,\cdot),f\rangle-\langle\mu_{1},f\rangle|\leq \epsilon/5
		\end{equation}\par 
		Combining (\ref{eq 5.1})-(\ref{eq 5.4}), we obtain that $|\langle\mu_1,f\rangle-\langle\mu_2,f\rangle|\leq\frac{4}{5}\epsilon,$ contrary to the definition of $\epsilon.$

	$\mathbf{Step\;3.}$ Finally, for any $x\in\text{supp }\mu$, we know that 	$\{Q_t(x,\cdot)\}_{t\in T}$ weakly converges to some invariant measure $\nu$. Assume that, contrary to our proposition, there exists some $x\in\text{supp }\mu$ such that $\nu\neq\mu$. Then we choose $f\in L_b(\mathcal{X})$ and $\epsilon>0$ such that $|\langle \mu,f\rangle-\langle \nu,f\rangle|>\epsilon$. By the Birkhoff's ergodic theorem, there exists $A\subset \mathcal{X},\;\mu(A)=1$, such that 
		\begin{equation}\label{eq 3.1}
		\lim\limits_{t\rightarrow\infty}Q_tf(y)=\langle \mu,f\rangle,\quad\text{for all }y\in A.
		\end{equation}\par 
		Note that $x\in\text{supp }\mu,$ we may find $y_k\in A$ such that $y_k\rightarrow x$ as $k\rightarrow\infty.$ Due to the eventual continuity of 	$\{P_t(x,\cdot)\}_{t\in T}$, there exists some $N\in\mathbb{N}$ sufficiently large such that for $n\geq N,$
		\begin{equation}\label{eq 3.2}
		\limsup\limits_{t\rightarrow\infty}|Q_tf(x)-Q_tf(y_n)|\leq\epsilon/2.
		\end{equation}\par 
		On the other hand,  $\{Q_t(x,\cdot)\}_{t\in T}$ weakly converges to $\nu,$ consequently,
		\begin{equation}\label{eq 3.3}
		\lim\limits_{t\rightarrow\infty}Q_tf(x)=\lim\limits_{t\rightarrow\infty}\langle Q_t(x,\cdot),f\rangle=\langle \nu,f\rangle.
		\end{equation}
		By (\ref{eq 3.1})-(\ref{eq 3.3}), it follows that $|\langle\mu,f\rangle-\langle \nu,f\rangle|\leq\epsilon/2,$ which contradicts with the definitions of $f$ and $\epsilon.$ This completes the proof.
	\end{prof}

    \section*{Acknowledgements} We would like to thank  Professor Fu-zhou Gong  and Professor Yuan Liu for their helpful comments and suggestions.

\end{document}